\documentclass{article}
\usepackage{amssymb,amsmath,amsthm}
\newtheorem{theorem}{Theorem}
\newtheorem{definition}{Definition}
\newtheorem{lemma}{Lemma}
\newtheorem{proposition}{Proposition}
\newtheorem{corollary}{Corollary}
\newtheorem{remark}{Remark}
\date{}
\numberwithin{equation}{section} \numberwithin{theorem}{section}
\numberwithin{lemma}{section} \numberwithin{corollary}{section}
\numberwithin{remark}{section} \numberwithin{proposition}{section}
\numberwithin{definition}{section}

\begin{document}
% let's define a new macro
\newcommand{\n}{\noindent}
\newcommand{\vs}{\vskip}

\title{Removable Sets for H\"{o}lder Continuous $p(x)$-Harmonic Functions }

\author{A. Lyaghfouri \\
Fields Institute, 222 College Street\\
Toronto M5T 3J1, Canada} \maketitle

\begin{abstract} We establish that a closed set $E$
is removable for  $C^{0,\alpha}$ H\"{o}lder continuous $p(x)$-harmonic functions in a
bounded open domain $\Omega$ of $\mathbb{R}^n$, $n\geq 2$, provided that
for each compact subset $K$ of $E$, the $(n-p_K+\alpha(p_K-1))$-Hausdorff measure of $K$
is zero, where $\displaystyle{p_K=\max_{x\in K} p(x)}$.
\end{abstract}

\n MSC: 35J60, 35J70 

\vs 0.2cm

\n Key words : $p(x)$-harmonic functions, Variable Exponents Sobolev Spaces,
H\"{o}lder continuity, Hausdorff Measure, Removable Sets.

\section{Introduction}\label{S:intro}

Let  $\Omega$ be an open bounded domain of $\mathbb{R}^n$,
$n\geqslant 2$, and let $p$ be a measurable real valued function defined in
$\Omega$ satisfying for some positive numbers $p_{-}$ and $p_{+}$
\begin{equation}\label{0.2}
1 < p_{-}  \leq p(x) \leq p_{+}<\infty \quad\text{a.e. } x\in \Omega.
\end{equation}

\vs0.3cm\n We first recall the definitions of Lebesgue  and Sobolev
spaces with variable exponents $L^{p(x)} (\Omega)$ and $W^{1,p(x)} (\Omega)$ (see for example \cite{[AS]},
\cite{[FZ]} and \cite{[KR]})

$$L^{p(x)} (\Omega) = \Big\{u : \Omega \rightarrow \mathbb{R} ~\mbox{ measurable}~ /~\rho(u)=\int_{\Omega}
|u(x)|^{p(x)} < + \infty~\Big\}$$
is equipped with the Luxembourg
norm $\displaystyle{\|u\|_{p(x)} = \inf\Big\{\lambda > 0~:~ \rho\left(\left|\frac{u(x)}{\lambda}\right|\right)\leqslant 1~\Big\}}.$

$$W^{1,p(x)} (\Omega) = \Big\{u \in L^{p(x)}\Omega)~ /~ \nabla u \in
\big(L^{p(x)}(\Omega)\big)^n\Big\}.$$

\n Setting $\|\nabla u\|_{p(x)} = \displaystyle{\sum^{n}_{i=1}
\left\|\frac{\partial u} {\partial x_{i}}\right\|_{p(x)}},$ for $u\in W^{1,p(x)}(\Omega)$, then $\|u\|_{1,p(x)} = \|u\|_{p(x)} + \|\nabla u\|_{p(x)}$ is a norm,
and $(W^{1,p(x)}(\Omega),\|.\|_{1,p(x)})$ is a separable and reflexive Banach space.

\n The space $W_0^{1,p(x)}(\Omega)$ is defined as the closure of
$C_0^{\infty}(\Omega)$ in $W^{1,p(x)} (\Omega)$.

\vs 0.3cm

\n If $\Omega$ has Lipschitz boundary and $p$ satisfies for some $L > 0$
\begin{equation}\label{1.3}
- |p(x) - p(y)| \log |x - y| \leqslant L \qquad \forall x, y \in
{\overline{\Omega}},\end{equation}

\n then $C^{\infty}(\overline{\Omega})$ is dense in $W^{1,p(x)}(\Omega)$ and
$W_0^{1,p(x)}(\Omega)=W^{1,p(x)}(\Omega)\cap W_0^{1,1}(\Omega)$ (see \cite{[FZ]}).

\vs 0.3cm\n We shall need three definitions:
\begin{definition}\label{d1.1}
We say that a function $u\in W_{loc}^{1,p(x)} (\Omega)$ is a $p(x)$-harmonic function
in an open subset ${\cal O}$ of $\Omega$ and write
$\Delta_{p(x)} u=0$ in ${\cal O}$, if it satisfies
$$\int_{\cal O}|\nabla u|^{p(x)-2} \nabla u . \nabla \zeta
  dx = 0
\qquad  \forall \zeta \in {\cal D}({\cal O}).$$
\end{definition}

\begin{definition}\label{d1.2}
Let $E$ be a closed subset of $\Omega$. We say that $E$ is a removable set for
$C^{0,\alpha}$ H\"{o}lder continuous $p(x)$-harmonic functions, if
for any $C^{0,\alpha}$ H\"{o}lder continuous function $u$ in $\Omega$:
$$\text{If} ~~u~~ \text{is } p(x)\text{-harmonic in }\Omega\setminus E, \text{ then } u \text{ is }
p(x)\text{-harmonic in }\Omega.$$
\end{definition}

\begin{definition}\label{d1.3}
Let $F$ be a subset of $\mathbb{R}^n$ and $s$ a positive real number. The $s$-Hausdorff measure of $F$, denoted by $H^s(F)$, is defined by
$$H^s(F)=\lim_{\delta\rightarrow 0}H_\delta^s(F)=\sup_{\delta>0}H_\delta^s(F),~\text{where for}~ \delta>0,$$
$$H_\delta^s(F)=\inf\Big\{~\sum_{j=1}^\infty\alpha(s)\Big({{diam (C_j)}\over
2}\Big)^s~|~ F\subset \displaystyle{\bigcup_{j=1}^\infty
}C_j,~diam (C_j)\leqslant \delta~\Big\}$$

\n and $\alpha(s)=\displaystyle{{\pi^{s/2}}\over {\Gamma(s/2+1)}}$,
$\Gamma(s)=\displaystyle{\int_0^\infty e^{-t}t^{s-1}dt}$, for
$s>0$ is the usual Gamma function.
\end{definition}

\vs 0.3cm In this paper, we are concerned with giving a sufficient
condition for a closed subset $E$ of $\Omega$ in order to be removable for $C^{0,\alpha}$
H\"{o}lder continuous $p(x)$-harmonic functions in $\Omega$.
We recall that for $p$ constant, Kilpel\"{a}inen and Zhong \cite{[KZ]} proved that a closed subset $E$ of $\Omega$ is
removable for $C^{0,\alpha}$ H\"{o}lder continuous $p$-harmonic functions if and only if
the $(n-p+\alpha(p-1))$-Hausdorff measure of $E$ is zero. Trudinger and Wang \cite{[TW]} proved
the sufficiency  of this condition. For $p=2$, this result is due to Carleson \cite{[C]}.
For a result in the framework of $A$-harmonic functions, i.e. functions satisfying
$div\Big({{a(|\nabla u|)}\over{|\nabla u|}}\nabla u\Big)=0$, with $\displaystyle{A(t)=\int_0^t a(s)ds}$ and
$a\in C^1( (0,+\infty))\cap C^0([0,+\infty))$, we refer to \cite{[CL]}.

\vs 0.3cm The main result of the paper is the following theorem:

\begin{theorem}\label{T1.1} Let $E\subset\Omega$ be a closed set. Assume that $u$ is a continuous
function in $\Omega$, $p(x)$-harmonic in $\Omega\setminus E$, and  such that for some $\alpha\in(0,1)$
$$ |u(x)- u(y)|\leqslant L|x-y|^\alpha\quad \forall y\in \Omega,~\forall x\in E.$$

\n If for each compact subset $K$ of $E$, the
$(n-p_K+\alpha(p_K-1))$-Hausdorff measure of $K$ is zero, where $p_K=\max_{x\in K} p(x)$,
then $u$ is $p(x)-$harmonic in $\Omega$.
\end{theorem}

\vs 0,2cm

\n An immediate consequence of Theorem 1.1 is the following corollary:

\begin{corollary}\label{c1.1}
A closed subset $E$ of $\Omega$ is a removable set for $C^{0,\alpha}$ H\"{o}lder continuous
$p(x)$-harmonic functions, if for each compact subset $K$ of $E$,
the $(n-p_K+\alpha(p_K-1))$-Hausdorff measure of $K$ is zero, where $p_K$ is as defined above.
\end{corollary}

\vs 0,2cm

\begin{remark}\label{r1.1} We recall that Carleson \cite{[C]}  proved Corollary 1.1 for the
Laplace operator. For $p-$Laplace like operators, it was established by Kilpel\"{a}inen and Zhong \cite{[KZ]},
and also by Trudinger and Wang \cite{[TW]}, under the assumption that $u$ has an ${\cal A}-$superharmonic
extension to $\Omega$. A partial result has been also obtained in \cite{[CL]} for $A$-harmonic functions.
 \end{remark}

%%%%%%%%%%%%%%%%%%%%%%%%%%%%%%%%%%%%%%%%%%%%%%%%%%%%
%%%%%%%%%%%%%%%%%%%%%%%%%%%%%%%%%%%%%%%%%%%%%%%%%%%%
%%%%%%%%%%%%%%%%%%%%%%%%%%%%%%%%%%%%%%%%%%%%%%%%%%%%
%%%%%%%%%%%%%%%%%%%%%%%%%%%%%%%%%%%%%%%%%%%%%%%%%%%%

\section{Proof of the Main Result}\label{2}

\n First, we introduce the following obstacle problem, where $D$ is a smooth subdomain of $\Omega$ and $\phi\in W^{1,p(x)}(D)$
$$P(\phi,D)
  \begin{cases}
   \text{Find}~ v\in \mathcal{F}=\{~\zeta\in W^{1,p(x)}(D)~/~\zeta\geqslant \phi~\text{in}~\Omega~\text{and}~\zeta-\phi\in W_0^{1,p(x)}(D)~ \}, & \\
    \displaystyle{\int_{D} |\nabla v|^{p(x)-2} \nabla v. \nabla (\zeta- v) dx \geqslant 0}  ~~\hbox{ for all } \zeta\in \mathcal{F}.&
  \end{cases}
$$
Then we have:

\begin{proposition}\label{p2.1} There exists a unique solution $v$ to the problem
$P(\phi,D)$. If $\phi$ is continuous in $D$, then so it is for $v$. Moreover
$-\Delta_{p(x)}v $ is a nonnegative measure and $v$ is $p(x)-$harmonic in $[v>\phi]$.
\end{proposition}

\n \emph{Proof.} The existence and uniqueness of a solution to $P(\phi,D)$ can easily be obtained
by standard techniques. For the rest, we refer to \cite{[HHKLM]} Theorem 10.
\qed

\vs 0,2cm

\n Next, we establish the following  lemma :

\begin{lemma}\label{l2.1}  Let $K$ be a compact subset of $\Omega$. Suppose that
$\phi\in W^{1,p(x)}(\Omega)$ is a continuous function such that
we have for some $L>0$
\begin{equation}\label{2.1}
\forall y\in \Omega~\forall x\in K~~ |\phi(x)- \phi(y)|\leqslant L |x-y|^\alpha.
\end{equation}
Let $v$ be the solution of the problem
$P(\phi,\Omega)$ and let $\mu= -\Delta_{p(x)}v $. Then there exit two positive constants $C_0$ and $R_0$ such that
\begin{equation}\label{2.2}
\mu( B_R(x)) \leqslant C_0 R^{n-p(x)+\alpha(p(x)-1)},\quad \forall R< R_1=\min\Big(R_0,{{dist(K, \partial\Omega)}\over 33}\Big),
~\forall x\in K.
\end{equation}
\end{lemma}

\n \emph{Proof.} Let $x_0\in K$ and $0<R< {{dist(K, \partial\Omega)}\over 33}$. We distinguish two cases:

\vs 0,2cm\n $\underline{1^{st} Case} :~~B_R(x_0)\cap[v=\phi]=\emptyset$.

\vs 0,2cm\n In this case, we have by Proposition 2.1  $~\mu(B_R(x_0))=0$.

\vs 0,2cm \n $\underline{2^{nd} Case} :~~\exists x_1\in B_R(x_0)\cap[v=\phi]$.

\vs 0,2cm\n In this case, we have $\mu(B_R(x_0))\leqslant \mu(B_{2R}(x_1))$ and it is enough to establish
(2.2) for $B_{2R}(x_1)$.

\n Note that since $\overline{B}_{33R}(x_0)\subset\subset\Omega$, we have
$\overline{B}_{32R}(x_1)\subset\subset\Omega$. Moreover if $v\in C^{0,\alpha}(\overline{B}_{8R}(x_1))$, then
we know (see \cite{[L]} Theorem 1.1) that there exist two positive constants $C_0$ and $r_0$ such that (2.2)
holds for $x=x_1$ provided that $2R<r_0$. We shall assume that $R< R_0=r_0/2$ and will show that
$\displaystyle{osc(v,\overline{B}_{8R}(x_1))=\max_{\overline{B}_{8R}(x_1)}v-\min_{\overline{B}_{8R}(x_1)}v}\leqslant CR^\alpha$
for some positive constant $C$.

\vs 0,2cm \n Let $\displaystyle{\omega_0=osc(\phi,\overline{B}_{32R}(x_1))=\max_{\overline{B}_{32R}(x_1)} \phi-\min_{\overline{B}_{32R}(x_1)} \phi}$.
Without loss of generality, we can assume that $v(x_1)=\phi(x_1)=0$, and we claim that $\omega_0 \leqslant L2^{\alpha+1}(1+32^\alpha)R^\alpha$.
Indeed since $\phi$ is continuous in $\Omega$, there exist $x_2,x_3\in \overline{B}_{32R}(x_1)$ such that
$\displaystyle{\phi(x_3)=\max_{\overline{B}_{32R}(x_1)} \phi}$ and $\displaystyle{\phi(x_2)=\min_{\overline{B}_{32R}(x_1)} \phi}$.
Then we have by the assumption (2.1) since $x_0\in K$
\begin{eqnarray}\label{2.3}
\omega_0&=& \phi(x_3)-\phi(x_2)=\phi(x_3)-\phi(x_0) +\phi(x_0)-\phi(x_2)\nonumber\\
&\leqslant& L |x_3-x_0|^\alpha+L |x_2-x_0|^\alpha\nonumber\\
&\leqslant& L (|x_3-x_1|+|x_1-x_0|)^\alpha+L(|x_2-x_1|+|x_1-x_0|)^\alpha\nonumber\\
&\leqslant& L 2^\alpha(|x_3-x_1|^\alpha+|x_1-x_0|^\alpha)+L2^\alpha(|x_2-x_1|^\alpha+|x_1-x_0|^\alpha)\nonumber\\
&\leqslant& L 2^\alpha((32R)^\alpha+R^\alpha+(32R)^\alpha+R^\alpha)=L2^{\alpha+1}(1+32^\alpha)R^\alpha.
\end{eqnarray}
Now $v+\omega_0$ is a nonnegative and $p(x)$-superharmonic function in $B_{32R}(x_1)$. Indeed
by Proposition 2.1, $-\Delta_{p(x)}(v+\omega_0)=-\Delta_{p(x)} v=\mu \geqslant 0$, and for $x\in B_{32R}(x_1)$
\begin{eqnarray*}
(v+\omega_0)(x)&=&v(x)+\max_{\overline{B}_{32R}(x_1)} \phi-\min_{\overline{B}_{32R}(x_1)} \phi\geqslant\phi(x)-\min_{\overline{B}_{32R}(x_1)} \phi +\max_{\overline{B}_{32R}(x_1)} \phi\nonumber\\
&\geqslant& \max_{\overline{B}_{32R}(x_1)} \phi\geqslant \phi(x_1)=0.
\end{eqnarray*}
Applying Lemma 6.4 of \cite{[AK]} to $v+\omega_0$, we get  for some constant $C_1>0$ and $0<q<{{n(p(x_1)-1)}\over{n-1}}$
\begin{eqnarray*}
\Bigg({1\over{|B_{16R}(x_1)|}}\int_{B_{16R}(x_0)} (v+\omega_0)^q dx\Bigg)^{1/q}\leqslant C_1\big(\inf_{B_{8R}(x_1)} (v+\omega_0)+R\big).
\end{eqnarray*}
Since by H\"{o}lder's inequality, we have
\begin{eqnarray*}
{1\over{|B_{16R}(x_1)|}}\int_{B_{16R}(x_1)} (v+\omega_0) dx\leqslant \Bigg({1\over{|B_{16R}(x_1)|}}\int_{B_{16R}(x_1)} (v+\omega_0)^q dx\Bigg)^{1/q} ,
\end{eqnarray*}
we get
\begin{eqnarray}\label{2.4}
{1\over{|B_{16R}(x_1)|}}\int_{B_{16R}(x_1)} (v+\omega_0) dx\leqslant C_1\big(\inf_{B_{8R}(x_1)} (v+\omega_0)+R\big).
\end{eqnarray}
\n Moreover $(v-\omega_0)^+$ is nonnegative and $p(x)$-subharmonic in $B_{32R}(x_1)$. Indeed
let $\zeta\in {\cal D}(B_{32R}(x_1))$, $\zeta \geqslant 0$, and $\epsilon>0$.
Using the fact that by Proposition 2.1 $\Delta_{p(x)} v=0$ in $[v>\phi]$,
and taking into account that $B_{32R}(x_1)\cap[v>\omega_0]\subset B_{32R}(x_1)\cap[v>\phi]$ because
$\displaystyle{-\min_{\overline{B}_{32R}(x_1)} \phi}\geqslant -\phi(x_1)=0$, we get
\begin{eqnarray*}
&&\int_{B_{32R}(x_1)} |\nabla v|^{p(x)-2} \nabla v . \nabla \big(\min\big(\zeta,{{(v-\omega_0)^+}\over\epsilon}\big)\big) dx\nonumber\\
&& =\int_{B_{32R}(x_1)\cap[v>\omega_0]} |\nabla v|^{p(x)-2} \nabla v . \nabla \big(\min\big(\zeta,{{(v-\omega_0)^+}\over\epsilon}\big)\big) dx= 0
\end{eqnarray*}
which can be written as
\begin{eqnarray*}
&&\int_{B_{32R}(x_1)\cap[\epsilon\zeta\leqslant(v-\omega_0)^+]} |\nabla (v-\omega_0)^+|^{p(x)-2} \nabla (v-\omega_0)^+ . \nabla \zeta dx\nonumber\\
&& =-{1\over \epsilon}\int_{B_{32R}(x_1)\cap[\epsilon\zeta>(v-\omega_0)^+]} |\nabla (v-\omega_0)^+|^{p(x)} dx\leqslant 0.
\end{eqnarray*}
Letting $\epsilon\rightarrow 0$, we obtain
\begin{eqnarray*}
\int_{B_{32R}(x_1)} |\nabla (v-\omega_0)^+|^{p(x)-2} \nabla (v-\omega_0)^+ . \nabla\zeta dx\leqslant 0
\end{eqnarray*}
which means that $\Delta_{p(x)}(v-\omega_0)^+\geqslant0$ in $B_{32R}(x_1)$.

\vs 0,2cm\n At this point, we remark that it is straightforward to adapt the proof of Lemma 6.6 of \cite{[AK]} to the function $(v-\omega_0)^+$, since
the proof uses only the fact that the function is $p(x)$-subharmonic. We obtain
for some constant $C_2>0$
\begin{eqnarray}\label{2.5}
\sup_{B_{8R}(x_1)} (v-\omega_0)^+&\leqslant &C_2\Big({1\over{|B_{16R}(x_1)|}}\int_{B_{16R}(x_1)} (v-\omega_0)^+ dx+R\Big)\nonumber\\
&\leqslant& C_2\Big({1\over{|B_{16R}(x_1)|}}\int_{B_{16R}(x_1)} (v+\omega_0)dx+R\Big).
\end{eqnarray}

\n Using (2.4)-(2.5), we get
\begin{eqnarray*}
\sup_{B_{8R}(x_1)} (v-\omega_0)^+&\leqslant &C_2\big(C_1\big(\inf_{B_{8R}(x_0)} (v+\omega_0)+R\big)+R\big)\leqslant C_3\big(\inf_{B_{8R}(x_0)} (v+\omega_0)+R\big)\nonumber\\
&\leqslant& C_3(\omega_0+R)  \qquad\text{since}~ \inf_{B_{8R}(x_1)} v \leqslant v(x_1)=0.
\end{eqnarray*}

\n We deduce that
\begin{equation*}
\sup_{B_{8R}(x_1)} v\leqslant \omega_0+\sup_{B_{8R}(x_1)} (v-\omega_0)\leqslant\omega_0+\sup_{B_{8R}(x_1)} (v-\omega_0)^+\leqslant (1+C_3)\omega_0+C_3 R
\end{equation*}
which leads by continuity to
\begin{equation}\label{2.6}
\max_{\overline{B}_{8R}(x_1)} v\leqslant (1+C_3)\omega_0+C_3.
\end{equation}

\n Now since $\phi(x_1)=0$, we have
\begin{eqnarray}\label{2.7}
-\min_{\overline{B}_{8R}(x_1)} v \leqslant -\min_{\overline{B}_{32R}(x_1)} \phi=\phi(x_1)-\min_{\overline{B}_{32R}(x_1)}\phi \leqslant \max_{\overline{B}_{32R}(x_1)}\phi-\min_{\overline{B}_{32R}(x_1)}\phi =\omega_0.
\end{eqnarray}

\n Since $\alpha\in(0,1)$, it follows from (2.3), (2.6) and (2.7) that
\begin{flushleft}
$~~\displaystyle{osc(v,B_{8R}(x_1))=\max_{\overline{B}_{8R}(x_1)} v-\min_{\overline{B}_{8R}(x_1)} v \leqslant(2+C_3)\omega_0+C_3 R}$

$~~\leqslant  L2^{\alpha+1}(1+32^\alpha)(2+C_3)R^\alpha+C_3 R=(L2^{\alpha+1}(1+32^\alpha)(2+C_3)+C_3R^{1-\alpha}) R^\alpha\leqslant C_4R^\alpha.$
\end{flushleft}

\qed

\vs 0,2cm \n \emph{Proof of Theorem 1.1}. Let $E$ be a closed subset of $\Omega$ such that
for each compact subset $K$ of $E$, the
$(n-p_K+\alpha(p_K-1))$-Hausdorff measure of $K$ is zero, where $p_K=\max_{x\in K} p(x)$.
Let $u$ be a continuous function in $\Omega$,
that is $p(x)$-harmonic in $\Omega\setminus E$, and  such that for all $y\in \Omega$ and $x\in E$
$$ |u(x)- u(y)|\leqslant L|x-y|^\alpha.$$

\n We would like to prove that $u$ is $p(x)$-harmonic in $\Omega$. Consider a smooth domain $D\subset\subset\Omega$.
By Proposition 2.1, there exists a unique continuous solution $v$ of the problem $P(u,D)$ such that $\mu=-\Delta_{p(x)}v$ is a
nonnegative Radon measure.
Let $K$ be a compact subset of $E\cap D$. Using Lemma 2.1, we have for a positive constant $C_0$
\begin{equation}\label{2.8}
\mu( B_R(x)) \leqslant C_0 R^{n-p(x)+\alpha(p(x)-1)}\qquad \forall x\in K,~\forall R<R_1=\min\Big(R_0,{{dist(K, \partial D)}\over 33}\Big).
\end{equation}
\n Let $\epsilon>0$. Since ${\cal H}^s(K)=0$ for $s=n-p_K+\alpha(p_K-1)$, there exists $\delta_0>0$ (see Definition 1.3) such that for all $\delta\in(0,\delta_0),~$
$0\leqslant H_\delta^s(K)\leqslant \epsilon.$ We deduce that for each $\delta\in(0,\delta_0)$, there exists
a family of sets $(C^\delta_j)$ such that $K\subset \displaystyle{\bigcup_{j=1}^\infty C_j^\delta}$,  $diam (C^\delta_j)\leqslant \delta$ and
\begin{equation}\label{2.9}
H_\delta^s(K)\leqslant\sum_{j=1}^\infty\alpha(s)\Big({{diam (C_j^\delta)}\over
2}\Big)^s<\epsilon.
\end{equation}
\n We assume naturally that for each $j$, $C^\delta_j\cap K\neq\emptyset$. So for each $j$, there exists an $x_j\in C^\delta_j\cap K$.  This leads to $C^\delta_j \subset B_{R_j}(x_j)$, with $R_j=diam(C^\delta_j)$. Obviously we can assume that for all $j$ $R_j<\min(1,R_1)$.
It follows from (2.8)-(2.9) that
\begin{eqnarray*}
\mu(K) &\leqslant& \sum_{j=1}^\infty \mu(C^\delta_j) \leqslant\sum_{j=1}^\infty \mu(B_{R_j}(x_j)) \leqslant C\sum_{j=1}^\infty R_j^{n-p(x_j)+\alpha(p(x_j)-1)}\nonumber\\
&\leqslant&C\sum_{j=1}^\infty R_j^{n-p_K+\alpha(p_K-1)}=C\sum_{j=1}^\infty (diam(C^\delta_j))^s\leqslant C{2^s\over{\alpha(s)}}\epsilon.
\end{eqnarray*}

\n Since $\epsilon$ is arbitrary, we get $\mu(K)=0$, which leads to $\mu(E\cap D)=0$.

\vs 0,2cm\n Next, we prove that $\mu(D\setminus E)=0$. Let $\zeta\in {\cal D}(D\setminus E)$, $\zeta \geqslant 0$,  $\epsilon>0$,
and set $\displaystyle{\zeta_\epsilon=\min\big(\zeta,{{v-u}\over\epsilon}\big)}$.

\n Given that we have $\mu=-\Delta_{p(x)}v=0$ in $[v>u]$ and $\zeta_\epsilon=0$
in $[v=u]$, we obtain

\begin{eqnarray}\label{e2.10}
\int_{D\setminus E} |\nabla v|^{p(x)-2}\nabla v. \nabla \zeta_\epsilon dx=\int_{(D\setminus E)\cap [v>u]} |\nabla v|^{p(x)-2}  \nabla v. \nabla \zeta_\epsilon dx=0.
\end{eqnarray}

\n Since $\Delta_{p(x)} u=0$ in $D\setminus E$, we have
\begin{eqnarray}\label{e2.10}
\int_{D\setminus E} |\nabla u|^{p(x)-2}  \nabla u. \nabla \zeta_\epsilon dx=0.
\end{eqnarray}

\n Subtracting (2.11) from (2.10), we get
\begin{eqnarray*}
\int_{D\setminus E} |\nabla v|^{p(x)-2} \nabla v-|\nabla u|^{p(x)-2} \nabla u\Big). \nabla \zeta_\epsilon dx=0
\end{eqnarray*}
which leads by the monotonicity of the vector function $|\xi|^{p(x)-2}\xi$ to
\begin{eqnarray*}
&&\int_{(D\setminus E)\cap[\epsilon\zeta\leqslant v-u]} \Big(|\nabla v|^{p(x)-2} \nabla v-|\nabla u|^{p(x)-2} \nabla u\Big). \nabla \zeta dx\nonumber\\
&& =-{1\over \epsilon}\int_{(D\setminus E)\cap[\epsilon\zeta>v-u]} \Big(|\nabla v|^{p(x)-2} \nabla v-|\nabla u|^{p(x)-2} \nabla u\Big). (\nabla v-\nabla u)dx\leqslant 0
\end{eqnarray*}

\n Letting $\epsilon\rightarrow 0$, we get
\begin{eqnarray*}
\int_{D\setminus E} \Big(|\nabla v|^{p(x)-2} \nabla v-|\nabla u|^{p(x)-2} \nabla u\Big). \nabla \zeta dx\leqslant 0.
\end{eqnarray*}

\n Now since $\Delta_{p(x)}u=0$ in $D\setminus E$, we obtain
\begin{eqnarray*}
\int_{D\setminus E} |\nabla v|^{p(x)-2}\nabla v. \nabla \zeta dx\leqslant 0.
\end{eqnarray*}

\n which means that $\mu=-\Delta_{p(x)}v\leqslant0$ in $D\setminus E$. We deduce that $\mu(D\setminus E)=0$,
and conclude that $\mu(D)=0$. Hence $\Delta_{p(x)}v=0$ in $D$.

\vs 0,2cm\n Similarly, we consider $w$ the solution of the obstacle problem $P(-u,D)$. In the same way we prove that
$\Delta_{p(x)}w=0$ in $D$. Now we have
$$\left\{
  \begin{array}{ll}
    \Delta_{p(x)}v=\Delta_{p(x)}(-w)=0, & \hbox{in}~~ D\\
    v=-w, & \hbox{on} ~~\partial D.
  \end{array}
\right.
$$

\n By the maximum principle, we get $v=-w$ in $D$. Since we have $-w\leqslant u\leqslant v$ in $D$, we obtain
 $u=v=-w$ in $D$ and $\Delta_{p(x)}u=0$ in $D$.
\qed

\vs 0,3cm
\section*{Acknowledgements} The author would like to thank the Fields Institute for
the facilities and excellent research conditions during his stay at this institution.

\end{document}